\documentclass[a4paper,11pt]{article}
\usepackage[top=30truemm,bottom=30truemm,left=30truemm,right=30truemm]{geometry}
\usepackage{amsmath}
\usepackage{amssymb}
\usepackage{float}
\usepackage{color}
\usepackage[dvipdfmx]{graphicx}
\usepackage{graphics}
\usepackage{graphicx}
\usepackage{tcolorbox}
\usepackage{latexsym}
\usepackage{ascmac}
\usepackage{framed}
\usepackage{fancybox}
\usepackage{theorem}
\usepackage{stmaryrd}
\usepackage{here}
\usepackage{mathtools}
\usepackage{appendix}

\mathtoolsset{showonlyrefs=true} 

\theoremstyle{break}

\makeatletter
\@addtoreset{equation}{section}

\makeatother

\makeatletter
\@addtoreset{figure}{section}

\makeatother

\title{
On global behavior of a some SIR epidemic model based on the Poincar\'e compactification
}

\author{
Yu Ichida
\thanks{Graduate School of Science and Technology, Meiji University, 1-1-1, Higashimita Tama-ku Kawasaki Kanagawa 214-8571, Japan, {\tt ichidayu@meiji.ac.jp} } $^{,}$ \footnote{Research Fellow of Japan Society for the Promotion of Science (JSPS Research Fellow)}
}

\begin{document}
\maketitle
\begin{abstract}
It is important to study the global behavior of solutions to systems of ordinary differential equations describing the transmission dynamics of infectious disease.
In this paper, we present a different approach from the Lyapunov function used in most of them.
This approach is based on the Poincar\'e compactification.
We then apply the method to a SIR endemic model as a test case, and discuss its effectiveness and the potential applications of this approach.
In addition, we refine the discussion of dynamics near the equilibrium, derive the asymptotic behavior, and mention its relation to the basic reproduction number.
\end{abstract}


\begin{center}
\end{center}

\section{Introduction}
\label{sec:sirC-int}
There have been many studies dealing with systems of ordinary differential equations (for short, ODEs) describing the transmission dynamics of infectious disease.
In particular, these fundamental questions are the global behavior of the solution for systems of ODEs, such as the number of infectious and their final state, and the derivation of the basic reproduction number, which plays an important role in understanding this behavior.
For details on the basic reproduction number and its derivation, for instance, see \cite{vand1} and references therein.
Let us denote the basic reproduction number by $\mathcal{R}_{0}$.

In many studies, the Lyapunov function is constructed and discussed in order to investigate the global behavior of the solution.
Note that understanding the global behavior of the solution is the same as studying the global stability of the equilibrium.
Global stability of an equilibrium is defined as the local asymptotic stability of the equilibrium and the convergence of all solutions of the system of ODE under consideration to the equilibrium.
The method of investigating the global stability of equilibria by means of the Lyapunov function is very powerful, but to the best of the author's knowledge, there is no theoretical system that constitutes the Lyapunov function in general (for instance, see \cite{Wiggins}).
In other words, when investigating the global behavior of solutions to complex systems of ODEs that incorporate more elements, constructing the Lyapunov function is itself a difficult problem. 
Therefore, investigating the global behavior of solutions is not an easy work in general.

In this paper, we present a method using the Poincar\'e compactification to study the global behavior of solutions.
We then apply the method to a mathematical model of infectious disease epidemics described by a system of ODE.
Note that this is essentially different from the analysis method using the Lyapunov function.
An overview of Poincar\'e compactification is given later in Section \ref{sec:sirC-PC}, which is one of the compactifications of the phase space (the embedding of $\mathbb{R}^{n}$ into the unit upper hemisphere of $\mathbb{R}^{n+1}$) (see, e.g. \cite{FAL, QTW, DNPE, RSSS, Matsue1, Matsue2}).
The most important feature of applying this method is that it allows us to investigate the global behavior of the system of ODEs of interest by revealing all its dynamics including infinity.
This method has been used, for instance, in the analysis of the Li\'enard equation (\cite{FAL} and references therein), in the classification of phase portraits of ODE systems derived from the Gray-Scotte model (\cite{TSJ21}), in the reconsideration of blow-up solutions of systems of ODEs in the view of dynamical systems (\cite{Matsue1, Matsue2}), and in the analysis of the behavior of characteristic and typical solutions of certain partial differential equations (\cite{QTW, DNPE, RSSS}).
To the best of the author's knowledge, there are no studies of its application to investigate the global behavior of solutions to ODEs related to infectious disease epidemics, such as the susceptible-infectious-recovered (SIR) model.
In this paper, we describe how to apply this method.
Using one of the most basic and representative mathematical models of infectious diseases, the SIR model such that it has birth and death terms (e.g., \cite{Hethcote}) as a test case, we show that the global behavior of the solution can be studied with this method.
It should be emphasized that this method can be used without constructing the Lyapunov function.

Furthermore, by refining the dynamics near the equilibria that are not located at infinity in this model (disease-free equilibria and epidemic equilibria, to be discussed later), we obtain the asymptotic behavior of the decay of infectious population for $\mathcal{R}_{0}<1$ and $\mathcal{R}_{0}=1$.
In particular, this is studied by applying the center manifold theory in the case that  $\mathcal{R}_{0}=1$.

This paper is organized as follows.
In the next section, we also briefly outline the Poincar\'e compactification that is the basis of this paper. 
In Section \ref{sec:sirC-model}, we describe the specific model described above and review the previously known methods for understanding the global behavior of solutions in that model.
This is based on the LaSalle's invariance principle with the Lyapunov function.
In Section \ref{sec:sirC-cal}, we apply the Poincar\'e compactification and concentrate on showing the effectiveness of this method for the concrete example presented in Section \ref{sec:sirC-model}.
Section \ref{sec:sirC-asy} shows the asymptotic behavior.
Finally, Section \ref{sec:sirC-con} is devoted to the conclusions and the possible applications of this method.

\section{Poincar\'e compactification}
\label{sec:sirC-PC}
In this section, we briefly introduce the Poincar\'e compactification. 
Here Section 2 of \cite{QTW, DNPE, RSSS} are reproduced.
Also, it should be noted that we refer \cite{FAL, Matsue1, Matsue2} for more details.

Let 
\[
X = P(S, I) \dfrac{\partial \mbox{}}{\partial S} 
+ Q(S, I) \dfrac{\partial \mbox{}}{\partial I}\]
be a polynomial vector field on $\mathbb{R}^{2}$, or in other words
\[
\begin{cases}
\dot{S} = P(S, I),
\\
\dot{I} = Q(S, I),
\end{cases}
\quad \left(\,\,\dot{}=\dfrac{d}{dt}\right),
\]
where $P$, $Q$ are polynomials of arbitrary degree in the variables $S$ and $I$.

First, we consider $\mathbb{R}^{2}$ as the plane in $\mathbb{R}^{3}$ defined by $(y_{1},y_{2},y_{3})=(S,I,1)$.
We consider the sphere 
\[
\mathbb{S}^{2} 
= \{ y \in \mathbb{R}^{3} \, |\, y_{1}^{2} + y_{2}^{2}+y_{3}^{2}=1\}
\]
which we call Poincar\'e sphere.
We divide the sphere into
\begin{align*}
 H_{+} &= \{ y \in \mathbb{S}^{2}\,|\,y_{3}>0\},\quad
 H_{-} = \{ y \in \mathbb{S}^{2}\,|\,y_{3}<0\},\quad
 \mathbb{S}^{1} = \{y \in \mathbb{S}^{2}\, | \, y_{3}=0\}.
\end{align*}

Let us consider the embedding of vector field $X$ from $\mathbb{R}^{2}$ to $\mathbb{S}^{2}$ given by $f^{+}:\mathbb{R}^{2} \to \mathbb{S}^{2}, \quad  f^{-}:\mathbb{R}^{2} \to \mathbb{S}^{2}$
where
\[
f^{\pm}(S, I):= \pm \left( \dfrac{S}{\Delta(S, I)},\dfrac{I}{\Delta(S,I)},\dfrac{1}{\Delta(S,I)} \right) 
\]
with $\Delta(S,I) := \sqrt{S^{2}+I^{2}+1}$.

Second, we consider six local charts on $\mathbb{S}^{2}$ given by $U_{k} = \{y \in \mathbb{S}^{2} \, | \, y_{k}>0\}$, $V_{k} = \{y \in \mathbb{S}^{2} \, | \, y_{k}<0\}$ for $k=1,2,3$.
Consider the local projection
\[
g^{+}_{k} : U_{k} \to \mathbb{R}^{2} \quad {\rm and} 
\quad g^{-}_{k} : V_{k} \to \mathbb{R}^{2} 
\]
defined as
\[
g^{+}_{k}(y_{1},y_{2},y_{3}) = - g^{-}_{k}(y_{1},y_{2},y_{3})
 = \left(\dfrac{y_{m}}{y_{k}},\dfrac{y_{n}}{y_{k}} \right) 
\]
for $m<n$ and $m,n \not= k$. 
The projected vector fields are obtained as the vector fields on the planes $\overline{U}_{k} = \{y \in \mathbb{R}^{3} \, | \, y_{k} = 1\}$ and $\overline{V}_{k} = \{y \in \mathbb{R}^{3} \, | \, y_{k} = -1\}$ for each local chart $U_{k}$ and $V_{k}$.
 We denote by $(x,\lambda)$ the value of $g^{\pm}_{k}(y)$ for any $k$.
 
\begin{figure}[t]
\centering
\includegraphics[width=7cm]{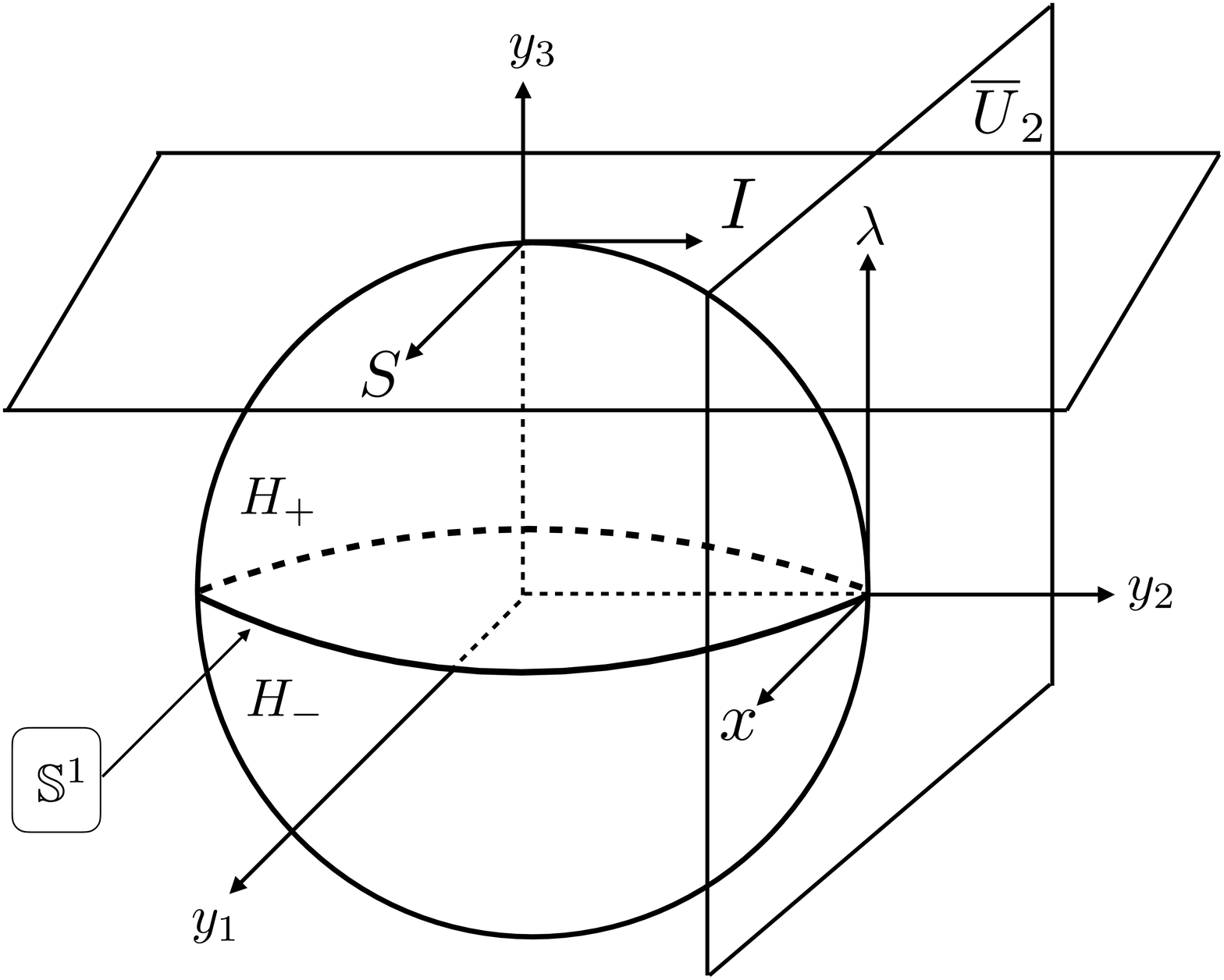}
\caption{Locations of the Poincar\'e sphere and chart $\overline{U}_{2}$}
\label{fig:ps-1}
\end{figure}

 For instance, it follows that
 \begin{align*}
 (x,\lambda) &:=(g^{+}_{2} \circ f^{+})(S,I) = \left ( \dfrac{S}{I},\dfrac{1}{I}\right) 
 \end{align*}
 therefore, we can obtain the dynamics on the local chart $\overline{U}_{2}$ 
 by the change of variables $S = x/\lambda$ and $I = 1/\lambda$.
 The locations of the Poincar\'e sphere, $(S, I)$-plane and $\overline{U}_{2}$ 
 are expressed as Figure \ref{fig:ps-1}. 
 The same is true for other local coordinates.
Now since we are considering non-negative solutions in the ODE system model, it is sufficient to consider $\overline{U}_{j}$ ($j=1,2$).

Throughout this paper, we follow the notations used here for the Poincar\'e compactification.
It is sufficient to consider the dynamics on $H_{+}\cup\mathbb{S}^{1}$, which is called Poincar\'e disk, to obtain our main results.

\section{The model as a test case}
\label{sec:sirC-model}
In this section, we present a model as a test case to verify the effectiveness of above the method based on the Poincar\'e compactification.
Before testing this method, we briefly review the global behavior obtained by using the LaSalle's invariance principle with the Lyapunov function that are known.

\subsection{The model}
\label{sub:sirC-mod1}
In this paper, the following model is treated as a concrete example (see \cite{Saitou2}):
\begin{equation}
\begin{cases}
\dot{S}=A-\beta S(t)I(t)-\mu S(t),
\\
\dot{I}=\beta S(t)I(t)-q I(t)-\mu I(t),
\\
\dot{R}= qI(t) -\mu R(t),
\end{cases}
\quad \left(\,\,\dot{}=\dfrac{d}{dt}\right).
\label{eq:sirC1}
\end{equation}
This model is a SIR ODE system model with birth and death terms, and is one of the classical and representative mathematical models for infectious disease.
To be more precise, we modify the Kermack-McKendrick model (for instance, see \cite{Murray}), which is the most typical SIR model, to introduce the fertility rate $A$ as inputs and the natural mortality rate $\mu$.
The model with $A=\mu$ is described, for instance, in \cite{Hethcote}.
The model \eqref{eq:sirC1} assumes that infection has no effect on births or deaths and that there is no vertical transmission.

$S(t)$ represents the susceptible population, $I(t)$ the infectious population, and $R(t)$ the recovered population, which are non-negative functions.
In this model, we impose initial conditions such as $S(0)=S_{0}>0$, $I(0)=I_{0}>0$, and $R(0)=R_{0}\ge 0$.
The constants $A, \beta, \mu, q$ are positive.
Furthermore, in \eqref{eq:sirC1}, the first and second equations do not contain any information about $R(t)$.
Therefore, it is sufficient to consider the following 2D: 
\begin{equation}
\begin{cases}
\dot{S}=A-\beta S(t)I(t)-\mu S(t),
\\
\dot{I}=\beta S(t)I(t)-q I(t)-\mu I(t),
\end{cases}
\quad \left(\,\,\dot{}=\dfrac{d}{dt}\right).
\label{eq:sirC2}
\end{equation}
The basic reproduction number for this model is known to be $\mathcal{R}_{0}=A\beta/\mu(q+\mu)$.

\subsection{Known results}
\label{sub:sirC-mod2}
\eqref{eq:sirC2} has the following equilibria:
\begin{align*}
E_{0}: (S, I)=\left( A\mu^{-1}, 0\right),
\quad
E_{*}: \left( (q+\mu)\beta^{-1}, \mu\beta^{-1}(\mathcal{R}_{0}-1)  \right)
=:(S_{*}, I_{*}),
\end{align*}
where $E_{0}$ is the disease-free equilibrium and $E_{*}$ is the endemic equilibrium.
Note that $E_{*}$ exists only if $\mathcal{R}_{0}>1$.
The Jacobian matrices of the vector filed \eqref{eq:sirC2} at these equilibria are
\[
E_{0}: \left(
\begin{array}{cc}
-\mu & -A\beta/\mu
\\
0 & [A\beta-\mu(q+\mu)]/\mu
\end{array}
\right),
\quad
E_{*}: \left(
\begin{array}{cc}
-A\beta/(q+\mu) & -(q+\mu)
\\
\mu(\mathcal{R}_{0}-1) & 0
\end{array}
\right).
\]
If $\mathcal{R}_{0}<1$, then, $E_{0}$ is a sink.
If $\mathcal{R}_{0}>1$, then $E_{0}$ is a saddle and $E_{*}$ is an asymptotically stable.

According to \cite{Saitou2}, if $\mathcal{R}_{0}>1$, then $E_{*}$ is a globally asymptotically stable, and if $\mathcal{R}_{0}\le 1$, then $E_{0}$ is a globally asymptotically stable.
When $\mathcal{R}_{0}>1$, we construct the following the Lyapunov function in region $\Omega$, which is shown by LaSalle's invariance principle:
\[
V(S,I)=S-S_{*}-S_{*}\log \{S(S_{*})^{-1}\}+I-I_{*}-I_{*}\log\{I(I_{*})^{-1}\}.
\]
Here, we set $\Omega=\{ (S,I) \mid S>0, \,\, I>0 \}$.
On the other hand, the case $\mathcal{R}_{0}\le 1$ can be shown similarly by using the following function:
\[
V(S,I)=S-A\mu^{-1}-A\mu^{-1}\log (A^{-1}\mu S)+1.
\]
Thus, the method using the Lyapunov function does not have a general theory of construction, and in this model we have to construct two functions.

\subsection{The case $\mathcal{R}_{0}=1$: Application of the center manifold theory}
\label{sub:sirC-mod3}
We focus on the Jacobian matrix of the vector field \eqref{eq:sirC2} at $E_{0}$ in the case that $\mathcal{R}_{0}=1$.
It has the real distinct eigenvalues $-\mu$ and $0$.
Therefore, the center manifold theory (e.g. \cite{carr}) is applicable to study the dynamics of \eqref{eq:sirC2}.
In \cite{Hethcote, Saitou2}, there is no discussion of this case.

We set 
\[
S(t)=A\mu^{-1}+U(t), \quad I(t)=0+V(t)
\]
and 
\[
\begin{cases}
\dot{U}=-\mu U-A\beta\mu^{-1}V-\beta UV,
\\
\dot{V}= \beta UV
\end{cases}
\]
hold.
Note that we are only shifting $E_{0}$ to the origin in this transformation.

The eigenvectors correspond to each eigenvalue are 
\[
{\mathbf{v}}_{1}=\left(\begin{array}{cc}
1 \\ 0
\end{array}\right),
\quad
{\mathbf{v}}_{2}=\left(\begin{array}{cc}
-A\beta \mu^{-1} \\ \mu
\end{array}\right).
\]
We set a matrix $T$ as $T=({\mathbf{v}}_{1}, {\mathbf{v}}_{2})$.
Then we obtain
\begin{align*}
\dfrac{d}{dt}\left(\begin{array}{cc}
U \\ V
\end{array}\right)
&=\left(\begin{array}{cc}
-\mu & -A\beta \mu^{-1}
\\
0 &  0
\end{array}\right)
\left(\begin{array}{cc}
U \\ V
\end{array}\right)
+
\left(\begin{array}{cc}
-\beta UV \\ \beta UV
\end{array}\right)
\\
&=T\left(\begin{array}{cc}
-\mu & 0
\\
0 &  0
\end{array}\right)T^{-1}
\left(\begin{array}{cc}
U \\ V
\end{array}\right)
+
\left(\begin{array}{cc}
-\beta UV \\ \beta UV
\end{array}\right).
\end{align*}
Let
\[
\left(\begin{array}{cc}
\tilde{U}
\\
\tilde{V}
\end{array}\right)
=T^{-1}\left(\begin{array}{cc}
U
\\
V
\end{array}\right).
\]
We then obtain the following system:
\[
\begin{cases}
\dot{\tilde{U}}= -\mu \tilde{U}+(A\beta^{2}\mu^{-1}-\beta \mu)\tilde{U}\tilde{V}+(A\beta^{2}-A^{2}\beta^{3}\mu^{-2})\tilde{V}^{2},
\\
\dot{\tilde{V}}=\beta \tilde{U}\tilde{V}-A\beta^{2}\mu^{-1}\tilde{V}^{2}.
\end{cases}
\]
The center manifold theory is applicable to study the dynamics of above system.
It implies that there exists a function $h(\tilde{V})$ satisfying 
\[
h(0)=\dfrac{dh}{d\tilde{V}}(0)=0
\]
such that the center manifold of the origin for above system is locally represented as $\{ (\tilde{U}, \tilde{V}) \mid \tilde{U}(t)=h(\tilde{V}(t)) \}$.
Differentiating it with respect to $t$, we have
\[
-\mu h+(A\beta^{2}\mu^{-1}-\beta \mu)h\tilde{V}+(A\beta^{2}-A^{2}\beta^{3}\mu^{-2})\tilde{V}^{2}
=\dfrac{dh}{d\tilde{V}} \left\{ \beta h\tilde{V}-A\beta^{2}\mu^{-1}\tilde{V}^{2} \right\}.
\]
Therefore, we obtain that the approximation of the (graph of) center manifold is 
\begin{equation}
\left\{ (S, I) \mid S(t)=A\mu^{-1}-A\beta\mu^{-2}I(t)+O(I^{2})  \right\}
\label{eq:sirC3}
\end{equation}
and the dynamics of \eqref{eq:sirC2} near $E_{0}$ is topologically equivalent to the dynamics of the following equation:
\begin{equation}
\dot{I}=-A\beta^{2}\mu^{-2}I^{2}+O(I^{3}).
\label{eq:sirC4}
\end{equation}
The above discussion has been made, for instance, in \cite{DNPE, DNLA}.

Here, the dynamics around $E_{0}$ is strictly different for $\mathcal{R}_{0}<1$ and $\mathcal{R}_{0}=1$.
Note, however, that in the later conclusions on global behavior and in Figure \ref{fig:sir-comp2}, they are treated together in $\mathcal{R}_{0}\le 1$ since they share the common feature of being attracted to $E_{0}$.

\section{Application of the Poincar\'e compactification}
\label{sec:sirC-cal}
In this model, we have a two-dimensional system, and we can use the Poincar\'e-Bendixson theorem (e.g. \cite{Wiggins}).
That is, if we know that there is no trajectory toward the equilibrium at infinity, we can show the global stability of the equilibrium without constructing the Lyapunov function, since it is attracted to the bounded equilibrium.

\subsection{Dynamics on the local charts}
\label{sub:sir-comp4-2}
First, to obtain the dynamics on the chart $\overline{U}_{2}$, we introduce the coordinates $(\lambda, x)$ by the formulas
\[
S(t)=x(t)/\lambda(t), \quad I(t)=1/\lambda(t).
\]
Then, we have
\[
\begin{cases}
\dot{\lambda}= -\beta x+(q+\mu) \lambda, 
\\
\dot{x}=A\lambda-\beta \lambda^{-1}x-\beta \lambda^{-1}x^{2}+qx,
\end{cases}
\quad \left(\,\,\dot{}=\dfrac{d}{dt}\right).
\]
By using the time-scale desingularization $d\tau/ dt=\lambda^{-1}$, we can obtain
\begin{equation}
\begin{cases}
\lambda_{\tau}=-\beta \lambda x+(q+\mu) \lambda^{2},
\\
x_{\tau}=A\lambda^{2}-\beta x-\beta x^{2}+q\lambda x,
\end{cases}
\label{eq:sirC5}
\end{equation}
where $\lambda_{\tau}=d\lambda /d\tau$ and $x_{\tau}=dx/d\tau$.
The equilibrium of the system \eqref{eq:sirC5} on $\{\lambda=0, x\ge 0\}$ is $E_{1}: (\lambda, x)=(0,0)$.
By calculating the Jacobian matrix, we apply the center manifold theory.
Therefore, we obtain that the approximation of the (graph of) center manifold is 
\begin{equation}
\left\{ (\lambda, x) \mid x(\tau)=A\beta^{-1}\lambda^{2}+O(\lambda^{3})  \right\}
\label{eq:sirC5-1}
\end{equation}
and the dynamics of \eqref{eq:sirC5} near $E_{1}$ is topologically equivalent to the dynamics of the following equation:
\begin{equation}
\lambda_{\tau}=(q+\mu)\lambda^{2}-A\lambda^{3}+O(\lambda^{4}).
\label{eq:sirC5-2}
\end{equation}
In conclusion, this argues that the trajectories will never go to $E_{1}$.

Second, we consider the dynamics on the chart $\overline{U}_{1}$.
From the change of coordinates $S(t)=1/\lambda(t)$, $I(t)=x(t)/\lambda(t)$, and time-rescaling $d\tau/ dt=\lambda^{-1}$, we obtain
\begin{equation}
\begin{cases}
\lambda_{\tau}=-A\lambda^{3}+\beta \lambda x+\mu \lambda^{2},
\\
x_{\tau}= \beta x-q\lambda x-A\lambda^{2}x+\beta x^{2}.
\end{cases}
\label{eq:sirC6}
\end{equation}
The equilibrium of the system \eqref{eq:sirC6} on $\{\lambda=0, x\ge 0\}$ is $E_{2}: (\lambda, x)=(0,0)$.
By calculating the Jacobian matrix, we apply the center manifold theory.
However, it is not shown explicitly as \eqref{eq:sirC5-1} and \eqref{eq:sirC5-2} since the center manifold is not unique.
As we can see from the nullcline, we conclude that the trajectories will never go to $E_{2}$.

\subsection{Dynamics on the Poincar\'e disk}
\label{sub:sir-comp4-4}
\begin{figure}[t]
\centering
\includegraphics[width=10cm]{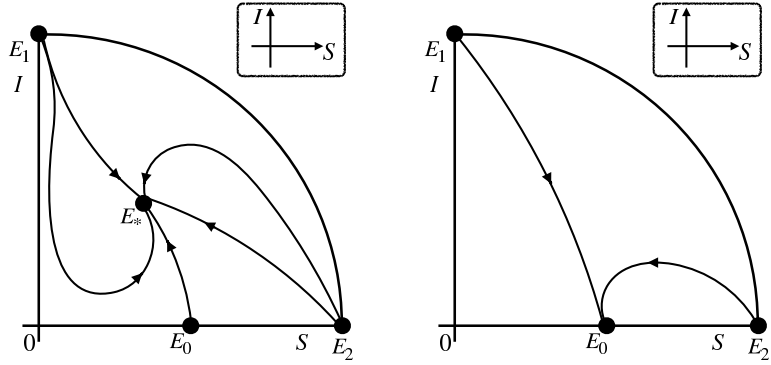}
\caption{
Schematic pictures of the dynamics on the Poincar\'e disk for \eqref{eq:sirC2}.
[Left: Case $\mathcal{R}_{0}>1$.] [Right: Case $\mathcal{R}_{0}\le 1$.]
Note that the circumference corresponds to $\{ \|(S, I)\|=+\infty \}$.
}
\label{fig:sir-comp2}
\end{figure}

Combining the dynamics on the charts $\overline{U}_{j}$ ($j=1,2$), we can obtain the dynamics on the Poincar\'e disk (see Figure \ref{fig:sir-comp2}).

We explain why the connected trajectories can be represented as shown in Figure \ref{fig:sir-comp2} in the both cases $\mathcal{R}_{0}>1$ and $\mathcal{R}_{0}\le 1$.
First, by computing $\dot{S}|_{S=0}$, $\dot{S}|_{I=0}$, $\dot{I}|_{S=0}$ and $\dot{I}|_{I=0}$, it is easy to see that both the Poincar\'e disks are an invariance set.
Second, bounded orbits cannot go to infinity by the local dynamics at infinity. 
Finally, the dynamics on these are determined from the Poincar\'e-Bendixson theorem.

Thus, we can understand the global behavior of the solution to \eqref{eq:sirC2} without constructing the Lyapunov function.
This means that the method introduced in Section \ref{sec:sirC-PC} is effective in \eqref{eq:sirC2}.

\section{Asymptotic behavior}
\label{sec:sirC-asy}
In this section, we derive the asymptotic behavior of $I(t)$ as $t\to\infty$ in the case that  $\mathcal{R}_{0}<1$ and $\mathcal{R}_{0}=1$, respectively.
We can also calculate the asymptotic behavior as $t\to +\infty$ in the case that  $\mathcal{R}_{0}>1$. 
See \cite{DNPE, RSSS} for similar argument.

First, from the Jacobian matrix at $E_{0}$ in the case that $\mathcal{R}_{0}<1$, the solution $I(t)$ can be approximated as follows
\begin{align*}
I(t)
&= C_{1}e^{[A\beta\mu^{-1}-(q+\mu)]t}(1+o(1))
 =C_{1}e^{(q+\mu)(\mathcal{R}_{0}-1)t}(1+o(1))
\quad {\rm{as}} \quad t\to +\infty. 
\end{align*}
That is, we obtain that $I(t)$ converges exponentially to $0$ with $t\to \infty$.
Here, $C_{1}$ is a positive constant.

Next, we examine the case that $\mathcal{R}_{0}=1$.
By considering terms up to the second order of \eqref{eq:sirC4} and solving for $I(t)$, we obtain 
\[
I(t) = (A\beta^{2}\mu^{-2}t +C_{2})^{-1}.
\]
This equation is the asymptotic behavior as $t\to \infty$.
Note that from $I(0)>0$, $C_{2}$ is a positive constant.

Although these results are same in that $I\to0$ is obtained as $t\to +\infty$, the decay is different.
Therefore, we can see that the smaller $\mathcal{R}_{0}$ is than $1$, the faster the convergence of $I(t)$ to $0$ as $t\to \infty$.
This is a result that gives a sense of the importance of making the basic reproduction number $\mathcal{R}_{0}$ small.

\section{Concluding remarks}
\label{sec:sirC-con}
The Poincar\'e compactification seems to be very convenient.
However, this method is not applicable in all cases.
To the best of the author's knowledge, local dynamics at infinity can be studied for phase spaces of dimension three or more, but it becomes difficult to immediately prove connecting orbits such as Subsection \ref{sub:sir-comp4-4} since the Poincar\'e-Bendixson theorem cannot be used.
Also, this method is not immediately applicable to non-polynomial vector fields in 2D system (for instance, see \cite{FAL}).
Since these problems are technical issues, their resolution will be future works.

\section*{Acknowledgments}
The author was partially supported by JSPS KAKENHI Grant Number JP21J20035.


\end{document}